\documentclass[12pt, twoside, leqno]{article}
\usepackage{amsmath,amsthm}
\usepackage{amssymb,latexsym}
\usepackage{enumerate}

\newtheorem{thm}{Theorem}[section]
\newtheorem{prop}[thm]{Proposition}
\newtheorem{cor}[thm]{Corollary}
\newtheorem{lem}[thm]{Lemma}

\theoremstyle{definition}

\newtheorem{rem}[thm]{Remark}

\numberwithin{equation}{section}

\usepackage{mathrsfs}
\usepackage{CJK}
\usepackage{amsmath}
\usepackage{fancyhdr}

\begin{document}

\baselineskip=17pt

\title{Polynomial maps with nilpotent Jacobians in dimension three II}
\author{ Dan Yan \footnote{ The author are supported by the Natural Science Foundation of Hunan Province (Grant No.2016JJ3085), the National Natural Science Foundation of China (Grant No.11601146) and the Construct Program of the Key Discipline in Hunan Province.}\\
Key Laboratory of HPCSIP,\\ College of Mathematics and Computer Science,\\
 Hunan Normal University, Changsha 410081, China \\
\emph{E-mail:} yan-dan-hi@163.com \\
}
\date{}

\maketitle

\renewcommand{\thefootnote}{}

\renewcommand{\thefootnote}{\arabic{footnote}}
\setcounter{footnote}{0}

%%%%%%%%
\begin{abstract} In the paper, we first classify all
polynomial maps of the form $H=(u(x,y),v(x,y,z), h(x,y))$ in the
case that $JH$ is nilpotent and $(\deg_yu,\deg_yh)\leq 3$, $H(0)=0$.
Then we classify all polynomial maps of the form
$H=(u(x,y,z),v(x,y,u), h(x,y))$ in the case that $JH$ is nilpotent
and $(\deg v(x,y,0),\allowbreak\deg h)\leq 3$, $H(0)=0$. Finally, we
classify polynomial maps of the form $H=(u(x,y,z),v(x,y,z), h(x,y))$
in certain conditions.

\end{abstract}
{\bf Keywords.} Jacobian Conjecture, Nilpotent Jacobian matrix, Polynomial maps\\
{\bf MSC(2010).} Primary 14E05;  Secondary 14A05;14R15 \vskip 2.5mm

\section{Introduction}

Throughout this paper, we will write ${\bf K}$ for algebraically
closed field and ${\bf K}[X]={\bf K}[x_1,x_2,\ldots,x_n]$ (${\bf
K}[\bar{X}]={\bf K}[x,y,z]$)for the polynomial algebra over ${\bf
K}$ with $n$ ($3$) indeterminates. Let
$F=(F_1,F_2,\ldots,F_n):{\bf{K}}^n\rightarrow{\bf{K}}^n$ be a
polynomial map, that is, $F_i\in{\bf{K}}[X]$ for all $1\leq i\leq
n$. Let $JF=(\frac{\partial F_i}{\partial x_j})_{n\times n}$ be the
Jacobian matrix of $F$.

The Jacobian Conjecture (JC) raised by O.H. Keller in 1939 in
\cite{1} states that a polynomial map
$F:{\bf{K}}^n\rightarrow{\bf{K}}^n$ is invertible if the Jacobian
determinant $\det JF$ is a nonzero constant. This conjecture has
been attacked by many people from various research fields, but it is
still open, even for $n\geq 2$. Only the case $n=1$ is obvious. For
more information about the wonderful 70-year history, see \cite{2},
\cite{3}, and the references therein.

In 1980, S.S.S.Wang (\cite{4}) showed that the JC holds for all
polynomial maps of degree 2 in all dimensions (up to an affine
transformation). The most powerful result is the reduction to degree
3, due to H.Bass, E.Connell and D.Wright (\cite{2}) in 1982 and
A.Yagzhev (\cite{5}) in 1980, which asserts that the JC is true if
the JC holds for all polynomial maps $X+H$, where $H$ is homogeneous
of degree 3. Thus, many authors study these maps and led to pose the
following problem.

 {\em (Homogeneous) dependence problem.} Let $H=(H_1,\ldots,H_n)\in
{\bf K}[X]$ be a (homogeneous) polynomial map of degree $d$ such
that $JH$ is nilpotent and $H(0)=0$. Whether $H_1,\ldots,H_n$ are
linearly dependent over ${\bf K}$?

The answer to the above problem is affirmative if rank$JH\leq 1$
(\cite{2}). In particular, this implies that the dependence problem
has an affirmative answer in the case $n=2$. D. Wright give an
affirmative answer when $H$ is homogeneous of degree 3 in the case
$n=3$ (\cite{6}) and the case $n=4$ is solved by Hubbers in
\cite{7}. M. de Bondt and A. van den Essen give an affirmative
answer to the above problem in the case $H$ is homogeneous and $n=3$
(\cite{8}). A. van den Essen finds the first counterexample in
dimension three for the dependence problem (\cite{9}). M. de Bondt
give a negative answer to the homogeneous dependence problem for
$d\geq 3$. In particular, he constructed counterexamples to the
problem for all dimensions $n\geq 5$ (\cite{10}). In \cite{18}, M.
Chamberland and A. van den Essen classify all polynomial maps of the
form $H=(u(x,y),v(x,y,z),h(u(x,y),v(x,y,z)))$ with $JH$ nilpotent.
In particular, they show that all maps of this form with $H(0)=0$,
$JH$ nilpotent and $H_1, H_2, H_3$ are linearly independent has the
same form as the counterexample that gave by A. van den Essen in
\cite{9} (up to a linear coordinate change). We classify all
polynomial maps of the form $H=(u(x,y),v(x,y,z),h(x,y,z))$ in the
case that $JH$ is nilpotent and $\deg_zv\leq 3$,
$(\deg_yu(x,y),\deg_yh(x,y,z))=1$ (\cite{13}) and classify all
polynomial maps of the form $H=(u(x,y),v(x,y,z),h(x,y))$ in the case
that $JH$ is nilpotent and $(\deg_yu(x,y),\deg_yh(x,y))\leq 2$ or
$\deg_yu(x,y)$ or $\deg_yh(x,y)$ is a prime number (\cite{14}).

In section 2, we classify all polynomial maps of the form
$H=(u(x,y),v(x,y,z),\allowbreak h(x,y))$ in the case that $JH$ is
nilpotent and $(\deg_yu,\deg_yh)\leq 3$, $H(0)=0$. Then, in section
3, we classify all polynomial maps of the form
$H=(u(x,y,z),v(x,y,u), \allowbreak h(x,y))$ in the case that $JH$ is
nilpotent and $(\deg v(x,y,0),\deg h)\leq 3$, $H(0)=0$. We prove
that $u,v,h$ are linearly dependent in the case that $JH$ is
nilpotent and $H$ has the form: $H=(u(x,y,z),v(x,y,z),h(x,y))$ with
$\deg_zv=1$ and $\deg_zu\geq 2$, $H(0)=0$ in section 4. The main
results in the paper are Theorem 2.2, Theorem 2.5, Theorem 3.2 and
Theorem 3.3, Theorem 4.1. We define $Q_{x_i}:=\frac{\partial
Q}{\partial x_i}$ and that $\deg_y f$ is the highest degree of $y$
in $f$.

\section{Polynomial maps of the form $H=(u(x,y),v(x,\allowbreak y,z), h(x,y))$}

In this section, we classify all polynomial maps of the form
$H=(u(x,y),v(x,y,z),\allowbreak h(x,y))$ in the case that $JH$ is
nilpotent and $(\deg_yu,\deg_yh)\leq 3$, $H(0)=0$.

\begin{lem}
Let $Q\in {\bf K}[x,y]$ and
$Q(x,y)=Q_r(x)y^r+Q_{r-1}(x)y^{r-1}+\cdots+Q_1(x)y+Q_0(x)$ with
$Q_r(x)\in {\bf K}^*$. If $Q_y|Q_x$ or $Q_y|Q+c$ for some $c\in {\bf
K}$, then $Q$ is a polynomial of $y+a(x)$ for some $a(x)\in {\bf
K}[x]$.
\end{lem}
\begin{proof}
Let $Q_r:=Q_r(x)$. Then we have
$Q_y=rQ_ry^{r-1}+(r-1)Q_{r-1}(x)y^{r-2}+\cdots+Q_1(x)$ and
$Q_x=Q_{r-1}'(x)y^{r-1}+\cdots+Q_1'(x)y+Q_0'(x)$. Thus, we have
$\deg_yQ_y\geq \deg_yQ_x$. We always view that the polynomials are
in ${\bf K}[x][y]$ with coefficients in ${\bf K}[x]$ when comparing
the coefficients of $y^i$.

Case I If $Q_y|Q_x$, then we have $$Q_x=b(x)Q_y~~~~~~~~~~~~~(2.1)$$
for some $b(x)\in {\bf K}[x]$. That is,
$$Q_{r-1}'(x)y^{r-1}+\cdots+Q_1'(x)y+Q_0'(x)=b(x)(rQ_ry^{r-1}+\cdots+Q_1(x))$$
Comparing the coefficients of $y^{r-1}$ of the above equation, we
have that $b(x)=\frac{Q_{r-1}'(x)}{rQ_r}$. Then let $\bar{x}=x$,
$\bar{y}=y+\frac{Q_{r-1}(x)}{rQ_r}$, it follows from equation
$(2.1)$ that $Q_{\bar{x}}=0$. That is, $Q(x,y)\in {\bf
K}[y+\frac{Q_{r-1}(x)}{rQ_r}]$. Let $a(x)=\frac{Q_{r-1}(x)}{rQ_r}$.
Then the conclusion follows.

Case II If $Q_y|(Q+c)$, then $Q+c=(e(x)y+d(x))Q_y$. That is,\\
$Q_ry^r+Q_{r-1}(x)y^{r-1}+\cdots+Q_1(x)y+Q_0(x)+c=(e(x)y+d(x))(rQ_ry^{r-1}+(r-1)Q_{r-1}(x)y^{r-2}+\cdots+Q_1(x))$~~~~~~~~~~~~~~~~~~~~~~~~~(2.2)\\
Comparing the coefficients of $y^r$ and $y^{r-1}$ of equation
$(2.2)$, we have that $e(x)=\frac{1}{r}$ and
$d(x)=\frac{Q_{r-1}(x)}{r^2Q_r}$. Then equation $(2.2)$ has the
following form: \\
$Q_ry^r+Q_{r-1}(x)y^{r-1}+\cdots+Q_1(x)y+Q_0(x)+c=(\frac{1}{r}y+\frac{Q_{r-1}(x)}{r^2Q_r})(rQ_ry^{r-1}+(r-1)Q_{r-1}(x)y^{r-2}+\cdots+Q_1(x))$~~~~~~~~~~~~~~~~~~~~~~~~~(2.3)\\

Claim: $Q_i(x)=C_r^i\frac{Q_{r-1}^{r-i}(x)}{r^{r-i}Q_r^{r-i-1}(x)}$
for $1\leq i\leq r-2$ and $C_r^i=\frac{r!}{i!(r-i)!}$.\\\\
Comparing the coefficients of $y^{r-2}$ of equation $(2.3)$, we have
the following equation:
$$Q_{r-2}(x)=C_r^{r-2}\frac{Q_{r-1}^2(x)}{r^2Q_r}$$
Suppose
$Q_{k+1}(x)=C_r^{k+1}\frac{Q_{r-1}^{r-k-1}(x)}{r^{r-k-1}Q_r^{r-k-2}(x)}$.
Then comparing the coefficients of $y^k$ of equation $(2.3)$, we
have the following equation:
$$Q_k(x)=\frac{1}{r}kQ_k(x)+\frac{Q_{r-1}(x)}{r^2Q_r}(k+1)Q_{k+1}(x).$$
That is,
$$Q_k(x)=C_r^k\frac{Q_{r-1}^r(x)}{r^{r-k}Q_r^{r-k-1}}.$$
Thus, we have that
$Q_i(x)=C_r^i\frac{Q_{r-1}^{r-i}(x)}{r^{r-i}Q_r^{r-i-1}(x)}$ for
$1\leq i\leq r-2$. Then equation $(2.3)$ has the following form:
$Q_0(x)+c=Q_1(x)\frac{Q_{r-1}(x)}{r^2Q_r}$. That is,
$$Q_0(x)=\frac{Q_{r-1}^r(x)}{r^rQ_r^{r-1}}-c.$$ Thus, we have
$Q(x,y)=Q_r(y^r+\frac{Q_{r-1}(x)}{Q_r}y^{r-1}+\cdots+\frac{Q_{r-1}^r(x)}{r^rQ_r^r})-c=Q_r(y+\frac{Q_{r-1}(x)}{rQ_r})^r-c$.
Therefore, we have $Q\in {\bf K}[y+\frac{Q_{r-1}(x)}{rQ_r}]$. Let
$a(x)=\frac{Q_{r-1}(x)}{rQ_r}$. Then the conclusion follows.
\end{proof}

\begin{thm}
Let $H=(u(x,y),v(x,y,z),h(x,y))$ be a polynomial map over ${\bf
K}[x,y,z]$. Assume that $H(0)=0$ and the components of $H$ are
linearly independent over ${\bf K}$. If $JH$ is nilpotent and
$(\deg_yu,\deg_yh)\leq 3$, then $u=g(ay+b(x))$,
$v=v_1z-a^{-1}b'(x)g(ay+b(x))-v_1l_2x$, $h=c_0u^2+l_2u$, where
$b(x)=v_1c_0ax^2+l_1x+\tilde{l}_2$; $v_1,c_0,a\in {\bf K}^*$;
$l_1,l_2,\tilde{l}_2\in {\bf K}$, $g(t)\in {\bf K}[t]$ and $g(0)=0$,
$\deg_tg(t)\geq 1$.
\end{thm}
\begin{proof}
If $(\deg_yu,\deg_yh)\leq 2$, then the conclusion follows from
Theorem 3.5 in \cite{14}. Let $v=v_dz^d+\cdots+v_1z+v_0$. Then it
follows from Lemma 3.2 in \cite{14} that $d=1$ and $v_1\in {\bf
K}^*$. Since $JH$ is nilpotent, so we have the following equations:
\begin{equation}
\nonumber
  \left\{ \begin{aligned}
  u_x+v_{0y} = 0~~~~~~~~~~~~~~~~~~~~~~~~(2.4) \\
  u_xv_{0y}-v_{0x}u_y-v_1h_y=0~~~~~~~~(2.5) \\
  v_1(u_xh_y-u_yh_x) = 0~~~~~~~~~~~~~~(2.6)
                          \end{aligned} \right.
  \end{equation}
It follows from equation $(2.4)$ that $u_x=-v_{0y}$. Thus, there
exists $P\in{\bf K}[x,y]$ such that
$$u=-P_y,~~v=P_x~~~~~~~~~~~~~~~~~~~~(2.7)$$
It follows from equation $(2.6)$ and Lemma 3.1 in \cite{14} that
there exists $q\in{\bf K}[x,y]$ such that
$$u,~h\in{\bf K}[q]~~~~~~~~~~~~~~~~~~~~~~~~~~~~~(2.8)$$
Since $\deg_yq|(\deg_yu,\deg_yh)$ and $(\deg_yu,\deg_yh)=3$, so we
have $\deg_yq=1$ or 3.

If $\deg_yq=1$, then the conclusion follows from the proof of
Theorem 2.8 in \cite{13}.

If $\deg_yq=3$, then $\deg_yq_y=2$. Let
$q(x,y)=q_3(x)y^3+q_2(x)y^2+q_1(x)y+q_0(x)$.\\
Claim: $q_3(x)\in{\bf K}^*$.\\
Let $P(x,y)=a_r(x)y^r+a_{r-1}(x)y^{r-1}+\cdots+a_1(x)y+a_0(x)$. It
follows from Lemma 3.4 in \cite{14} that $a_r(x)\in {\bf K}^*$. It
follows from equations $(2.7)$ and $(2.8)$ that
$u(x,y)=u(q)=-P_y=-(ra_r(x)y^{r-1}+\cdots+a_1(x))$. Thus, we have
$q_3(x)\in{\bf K}^*$. Substituting equations $(2.7)$ and $(2.8)$ to
equation $(2.5)$, we have the following equation:
$$q_y[v_1h'(q)+v_{0x}u'(q)]=-(u'(q)q_x)^2~~~~~~~~~~~~~(2.9)$$

(1) If $q_y$ is irreducible, then we have $q_y|q_x$ or $q_y|u'(q)$.
Since $u'(q)$ is a polynomial of $q$, so we have
$u'(q)=c_0(q+c_1)(q+c_2)\cdots(q+c_k)$ with $c_0\in{\bf K}^*$ and
$c_i\in{\bf K}$, $1\leq i\leq k$. Since $q_y$ is irreducible, so
there exists $i_0\in\{1,2,\ldots,k\}$ such that $q_y|q+c_{i_0}$ in
the case that $q_y|u'(q)$. It follows from Lemma 2.1 that $q$ is a
polynomial of $y+a(x)$ for some $a(x)\in{\bf K}[x]$. That is,
$u,h\in{\bf K}[y+a(x)]$. Then the conclusion follows from the proof
of Theorem 2.8 in \cite{13}.

(2) If $q_y$ is reducible, then
$$q_y=3q_3y^2+2q_2(x)y+q_1(x)=3q_3(y+b(x))(y+e(x))~~~~~~~~~~~(2.10)$$
for some $b(x),~e(x)\in{\bf K}[x]$. Clearly, $y+b(x)$, $y+e(x)$ are
irreducible. It follows from equation $(2.9)$ that we have the
following cases:\\
Case I $q_y|q_x$\\
Case II $q_y|u'(q)$\\
Case III $(y+b(x))|q_x$ and $(y+e(x))|u'(q)$\\
Case IV $(y+b(x))|q_x$ in the case $b(x)=e(x)$.\\
Case V $(y+b(x))|u'(q)$ in the case $b(x)=e(x)$.

Case I If $q_y|q_x$, then the conclusion follows from Lemma 2.1.

Case II If $q_y|u'(q)$, then $(y+b(x))(y+e(x))|u'(q)$. Since $u'(q)$
is a polynomial of $q$, so $u'(q)=c_0(q+c_1)(q+c_2)\cdots(q+c_k)$
for $c_0\in{\bf K}^*$ and $c_i\in{\bf K}$, $1\leq i\leq k$. Thus,
there exist $i_0,~j_0\in\{1,2,\ldots,k\}$ such that
$(y+b(x))|q+c_{i_0}$ and $(y+e(x))|q+c_{j_0}$. That is,
$$q_3y^3+q_2(x)y^2+q_1(x)y+q_0(x)+c_{i_0}=(y+b(x))(q_3y^2+d_1(x)y+d_0(x))~~~~~(2.11)$$
and
$$q_3y^3+q_2(x)y^2+q_1(x)y+q_0(x)+c_{j_0}=(y+e(x))(q_3y^2+f_1(x)y+f_0(x))~~~~~(2.12)$$
We always view that the polynomials are in ${\bf K}[x][y]$ with
coefficients in ${\bf K}[x]$ in the following arguments. Comparing
the coefficients of $y^2$, $y$ of equations $(2.11)$ and $(2.12)$,
we have the following equations:
$$q_2(x)=d_1(x)+q_3b(x)=f_1(x)+q_3e(x)~~~~~~~~~~~(2.13)$$
and
$$q_1(x)=d_0(x)+d_1(x)b(x)=f_0(x)+e(x)f_1(x)~~~~~(2.14)$$
It follows from equations $(2.13)$ and $(2.14)$ that
$$d_0(x)=q_1(x)-b(x)q_2(x)+q_3b^2(x)~~~~~~~~~~~~~~(2.15)$$
and
$$f_0(x)=q_1(x)-e(x)q_2(x)+q_3e^2(x)~~~~~~~~~~~~~~~(2.16)$$
Then equations $(2.11)$ and $(2.12)$ have the following forms:
$$q_0(x)+c_{i_0}=b(x)d_0(x)=b(x)(q_1(x)-b(x)q_2(x)+q_3b^2(x))~~~~~~~~(2.17)$$
and
$$q_0(x)+c_{j_0}=e(x)f_0(x)=e(x)(q_1(x)-e(x)q_2(x)+q_3e^2(x))~~~~~~~~~(2.18)$$
Then equation $(2.17)-(2.18)$ has the following form:
$$c_{i_0}-c_{j_0}=(b(x)-e(x))[q_1(x)-q_2(x)(b(x)+e(x))+q_3(b^2(x)+b(x)e(x)+e^2(x))].$$
Thus, we have $b(x)-e(x)\in {\bf K}$. That is, $e(x)=b(x)+c$ for
some $c\in {\bf K}$. Therefore, we have
$q_y=3q_3(y+b(x))(y+b(x)+c)=3q_3(y+b(x))^2+3q_3c(y+b(x))$. Then
$q(x,y)=q_3(y+b(x))^3+\frac{3}{2}q_3c(y+b(x))^2+f(x)$ for some
$f(x)\in {\bf K}[x]$. Since $(y+b(x))|(q+c_{i_0})$, so we have
$(y+b(x))|(f(x)+c_{i_0})$. That is, $f(x)=-c_{i_0}$. Thus, $q$ is a
polynomial of $y+b(x)$. Therefore, $u,h\in {\bf K}[y+b(x)]$. Then
the conclusion follows from the proof of Theorem 2.8 in \cite{13}\\

Case III Since $u'(q)=c_0(q+c_1)(q+c_2)\cdots(q+c_k)$ for $c_0\in
{\bf K}^*$ and $c_i\in{\bf K}$, $1\leq i\leq k$. Clearly, $y+e(x)$
is irreducible, so there exists $i_0\in \{1,2,\ldots,k\}$ such that
$(y+e(x))|(q+c_{i_0})$. That is,
$$q'_2(x)y^2+q_1'(x)y+q'_0(x)=(y+b(x))(q'_2(x)y+d(x))~~~~~~~~~~~~~~(2.19)$$
and
$$q_3y^3+q_2(x)y^2+q_1(x)y+q_0(x)+c_{i_0}=(y+e(x))(q_3y^2+w_1(x)y+w_0(x))~~~~~~(2.20)$$
where $d(x),w_1(x),w_0(x)\in{\bf K}[x]$. Comparing the coefficients
of $y$, $y^0$ of equation $(2.19)$, we have the following equations:
$$q'_1(x)-b(x)q'_2(x)=d(x)$$
and
$$q'_0(x)=b(x)d(x)=b(x)(q'_1(x)-b(x)q'_2(x))~~~~~~~~~(2.21)$$
Comparing the coefficients of $y$, $y^0$ of equation $(2.10)$, we
have the following equations:
$$\frac{2q_2(x)}{3q_3}=b(x)+e(x)~~~~~~~~~~~~~~~(2.22)$$
and
$$\frac{q_1(x)}{3q_3}=b(x)e(x)~~~~~~~~~~~~~~~~~(2.23)$$
Substituting equations $(2.22)$ and $(2.23)$ to equation $(2.21)$,
we have the following equation:
$$q'_0(x)=q_3(3b(x)b'(x)e(x)+\frac{3}{2}b^2(x)e'(x)-\frac{3}{2}b^2(x)b'(x))$$
Integrating the two sides of the above equation with respect to $x$,
we have the following equation:
$$q_0(x)=q_3(\frac{3}{2}b^2(x)e(x)-\frac{1}{2}b^3(x))+\bar{c}~~~~~~~~~~~~~(2.24)$$
for some $\bar{c}\in {\bf K}$. Comparing the coefficients of $y^2$,
$y$ and $y^0$ of equation $(2.20)$, we have the following equations:
\begin{equation}
\nonumber
  \left\{ \begin{aligned}
  q_2(x) = w_1(x)+q_3e(x)~~~~~~~~~~~~(2.25) \\
  q_1(x) = w_0(x)+w_1(x)e(x)~~~~~~~(2.26) \\
  q_0(x)+c_{i_0} = e(x)w_0(x)~~~~~~~~~~~(2.27)\\
                          \end{aligned} \right.
  \end{equation}
It follows from equations $(2.25)$ and $(2.26)$ that
$$w_0(x)=q_1(x)-e(x)q_2(x)+q_3e^2(x)~~~~~~~~~~~(2.28)$$
Substituting equations $(2.28)$, $(2.22)$ and $(2.23)$ to equation
$(2.27)$, we have the following equation:
$$q_0(x)=q_3(\frac{3}{2}b(x)e^2(x)-\frac{1}{2}e^3(x))-c_{i_0}~~~~~~~~~~~(2.29)$$
Comparing equation $(2.29)$ with equation $(2.24)$, we have
$\frac{1}{2}q_3[b^3(x)-3b^2(x)e(x)+3b(x)e^2(x)-e^3(x)]=\bar{c}+c_{i_0}$.
That is, $$\frac{1}{2}q_3(b(x)-e(x))^3=\bar{c}+c_{i_0}.$$ Thus, we
have $b(x)-e(x)\in {\bf K}$. That is, $b(x)=e(x)+\bar{\bar{c}}$ for
some $\bar{\bar{c}}\in {\bf K}$. It follows from equation $(2.10)$
that $q_y=3q_3(y+e(x))^2+3q_3\bar{\bar{c}}(y+e(x))$. Then we have
$q=q_3(y+e(x))^3+\frac{3}{2}q_3\bar{\bar{c}}(y+e(x))^2+\bar{f}(x)$
for some $\bar{f}(x)\in{\bf K}[x]$. Since $(y+e(x))|(q+c_{i_0})$, so
we have $(y+e(x))|(\bar{f}(x)+c_{i_0})$. Thus, we have
$\bar{f}(x)=-c_{i_0}$. Therefore, $q$ is a polynomial of $y+e(x)$.
That is, $u,h,q\in {\bf K}[y+e(x)]$. Then the conclusion follows
from the proof of Theorem 2.8 in \cite{13}.\\

Case IV, V If $b(x)=e(x)$, then it follows from equation $(2.10)$
that $q_y=3q_3(y+b(x))^2$. Thus, we have
$q=q_3(y+b(x))^3+\bar{\bar{f}}(x)$ for some $\bar{\bar{f}}(x)\in
{\bf K}[x]$.

If $(y+b(x))|q_x$, then
$(y+b(x))|(3q_3b'(x)(y+b(x))^2+\bar{\bar{f}}'(x))$. Thus, we have
$\bar{\bar{f}}'(x)=0$. That is, $\bar{\bar{f}}(x)\in {\bf K}$.
Therefore, we have $q\in {\bf K}[y+b(x)]$.

If $(y+b(x))|u'(q)$, then it follows from the arguments of Case II
that $(y+b(x))|q+c_{i_0}$ for some $c_{i_0}\in {\bf K}$. That is,
$(y+b(x))|(q_3(y+b(x))^3+\bar{\bar{f}}(x)+c_{i_0})$. Thus, we have
$\bar{\bar{f}}(x)=-c_{i_0}\in {\bf K}$. Therefore, $q\in{\bf
K}[y+b(x)]$.\\
Thus, $u,~h$ are polynomials of $y+b(x)$ in the two cases. Then the
conclusion follows from the proof of Theorem 2.8 in \cite{13}.
\end{proof}

\begin{rem}
We can replace the condition that $(\deg_yu(x,y),m)=1$ by the
condition $(\deg_yu(x,y),m)\leq 3$ in Theorem 2.10 and replace the
condition that $(m,n)=1$ by the condition $(m,n)\leq 3$ in Theorem
3.2 and Theorem 3.4 in \cite{13}.
\end{rem}

\begin{cor}
Let $H=(u(x,y),v(x,y,z),h(x,y))$ be a polynomial map over ${\bf
K}[x,y,z]$. Assume that $H(0)=0$ and the components of $H$ are
linearly independent over ${\bf K}$. If $JH$ is nilpotent and
$\deg_yu\leq 7$ or $\deg_yu=9$ or $\deg_yh\leq 7$ or $\deg_yh=9$,
then $H$ has the form of Theorem 2.2.
\end{cor}
\begin{proof}
Let $v=v_dz^d+\cdots+v_1z+v_0$. Then it follows from Lemma 3.2 in
\cite{14} that $d=1$ and $v_1\in {\bf K}^*$. Since $JH$ is
nilpotent, so we have the following equations:
\begin{equation}
\nonumber
  \left\{ \begin{aligned}
  u_x+v_{0y} = 0~~~~~~~~~~~~~~~~~~~~~~~~(2.4) \\
  u_xv_{0y}-v_{0x}u_y-v_1h_y=0~~~~~~~~(2.5) \\
  v_1(u_xh_y-u_yh_x) = 0~~~~~~~~~~~~~~(2.6)
                          \end{aligned} \right.
  \end{equation}
It follows from equation $(2.6)$ and Lemma 3.1 in \cite{14} that
there exists $q\in{\bf K}[x,y]$ such that
$$u,~h\in{\bf K}[q]~~~~~~~~~~~~~~~~~~~~~~~~~~~~~(2.8)$$
Thus, we have $\deg_yq|(\deg_yu,\deg_yh)$.

If $\deg_yu\leq 4$ or $\deg_yh\leq 4$, then the conclusion follows
from Corollary 3.9 in \cite{14}.

If $\deg_yu=5$ or $\deg_yu=7$ or $\deg_yh=5$ or $\deg_yh=7$, then
the conclusion follows from Theorem 3.3 in \cite{14}.

If $\deg_yu=6$ or 9 or $\deg_yh=6$ or 9, then $\deg_yq=1$ or 2 or 3
or 6 or 9.

Case I If $\deg_yq=1$ or 2 or 3, then the conclusion follows from
the proof Theorem 2.2.

Case II If $\deg_yq=6$, then $\deg_yu=6$ or $\deg_yh=6$.

(1) If $\deg_yq=\deg_yu$, then it follows from equation $(2.8)$ that
$u(x,y)=u(q)=\lambda q+\lambda_0$ with $\lambda\in{\bf K}^*$,
$\lambda_0\in{\bf K}$. That is,
$q=\lambda^{-1}u-\lambda^{-1}\lambda_0$. Thus, $h$ is a polynomial
of $u$. Then the conclusion follows from Theorem 2.1 in \cite{18}.

(2) If $\deg_yq=\deg_yh$, then it follows from the arguments of (1)
that $u$ is a polynomial of $h$. It follows from Corollary 2.3 in
\cite{14} that $u,v,h$ are linearly dependent. This is a
contradiction!

Case III If $\deg_yq=9$, then $\deg_yu=9$ or $\deg_yh=9$. That is,
$\deg_yq=\deg_yu$ or $\deg_yq=\deg_yh$. Then the conclusion follows
from the arguments of Case II.
\end{proof}

\begin{thm}
Let $H=(u(x,y),v(x,y,z),h(x,y))$ be a polynomial map over ${\bf
K}[x,y,z]$. Assume that $H(0)=0$ and the components of $H$ are
linearly independent over ${\bf K}$. If $JH$ is nilpotent and the
variety $V(Q_y)$ is irreducible for any $Q_y|(u_y,h_y)$ and
$\deg_yQ|(\deg_yu,\deg_yh)$, then $H$ has the form of Theorem 2.2.
\end{thm}
\begin{proof}
Let $v=v_dz^d+\cdots+v_1z+v_0$. Then it follows from Lemma 3.2 in
\cite{14} that $d=1$ and $v_1\in {\bf K}^*$. Since $JH$ is
nilpotent, so we have the following equations:
\begin{equation}
\nonumber
  \left\{ \begin{aligned}
  u_x+v_{0y} = 0~~~~~~~~~~~~~~~~~~~~~~~~(2.4) \\
  u_xv_{0y}-v_{0x}u_y-v_1h_y=0~~~~~~~~(2.5) \\
  v_1(u_xh_y-u_yh_x) = 0~~~~~~~~~~~~~~(2.6)
                          \end{aligned} \right.
  \end{equation}
It follows from equation $(2.4)$ that $u_x=-v_{0y}$. Thus, there
exists $P\in{\bf K}[x,y]$ such that
$$u=-P_y,~~v=P_x~~~~~~~~~~~~~~~~~~~~(2.7)$$
It follows from equation $(2.6)$ and Lemma 3.1 in \cite{14} that
there exists $q\in{\bf K}[x,y]$ such that
$$u,~h\in{\bf K}[q]~~~~~~~~~~~~~~~~~~~~~~~~~~~~~(2.8)$$
If $\deg_yq\leq 1$, then the conclusion follows from the proof of
Theorem 2.8 in \cite{13}.\\

Suppose $\deg_yq\geq 2$. Let
$q(x,y)=q_t(x)y^t+q_{t-1}(x)y^{t-1}+\cdots+q_1(x)y+q_0(x)$ with
$t\geq 2$. It follows from equation $(2.8)$ and Lemma 3.4 in
\cite{14} that $$q_t(x)\in{\bf
K}^*~~~~~~~~~~~~~~~~~~~~~~~~~~~~(2.30)$$ It follows from equations
$(2.5)$, $(2.7)$ and $(2.8)$ that
$$q_y[v_1h'(q)+v_{0x}u'(q)]=-(u'(q)q_x)^2~~~~~~~~~~~~~(2.31)$$
Since $u_y=u'(q)q_y$, $h_y=h'(q)q_y$, so we have $q_y|(u_y,h_y)$ and
$\deg_yq|(\deg_yu,\deg_yh)$. Thus, the variety $V(q_y)$ is
irreducible. That is, there exists an irreducible polynomial
$f(x,y)\in{\bf K}[x,y]$ such that
$$q_y=f^s(x,y)~~~~~~~~~~~~~~~~~~~~(2.32)$$
for some $s\in{\bf N}^*$. It follows from equation $(2.31)$ that
$q_y|(u'(q))^2q^2_x$. That is, $f^s(x,y)|(u'(q))^2q^2_x$. Since
$f(x,y)$ is irreducible, so we have $f(x,y)|u'(q)$ or
$f(x,y)|q_x$.\\

Case I Suppose $f(x,y)|u'(q)$. Since $u'(q)$ is a polynomial of $q$,
so it follows from the Fundamental Theorem of Algebra that
$u'(q)=c_0(q+c_1)(q+c_2)\cdots(q+c_k)$ for $c_0\in{\bf K}^*$ and
$c_i\in{\bf K}$, $1\leq i\leq k$. If $f(x,y)|u'(q)$, then there
exists $i_0\in\{1,2,\ldots,k\}$ such that $f(x,y)|(q+c_{i_0})$. That
is, $$q+c_{i_0}=f(x,y)M^{(1)}(x,y)~~~~~~~~~~~~~(2.33)$$ for some
$M^{(1)}(x,y)\in{\bf K}[x,y]$. Thus, it follows from equations
$(2.33)$ and $(2.32)$ that $q_y=f_yM^{(1)}+fM^{(1)}_y=f^s$. That is,
$$f_yM^{(1)}=f(f^{s-1}-M^{(1)}_y).$$
Since $f$ is irreducible and $f_y\neq 0$, so it follows from the above equation that
$f|M^{(1)}$. That is,
$$M^{(1)}(x,y)=f(x,y)M^{(2)}(x,y)~~~~~~~~~~~~(2.34)$$
for some $M^{(2)}(x,y)\in{\bf K}[x,y]$. Thus, it follows from
equations $(2.33)$ and $(2.34)$ that
$$q+c_{i_0}=f^2(x,y)M^{(2)}(x,y)~~~~~~~~~~~~~~~~(2.35)$$
Then it follows from equation $(2.35)$ that
$q_y=2ff_yM^{(2)}+f^2M^{(2)}_y=f^s$. That is,
$$2f_yM^{(2)}=f(f^{s-2}-M^{(2)}_y).$$
Since $f$ is irreducible and $f_y\neq 0$, so we have $f|M^{(2)}$. We can do step by
step until
$$q+c_{i_0}=f^{s+1}M^{(s+1)}$$
for some $M^{(s+1)}\in{\bf K}[x,y]$. Then it follows from the above
equation and equation $(2.32)$ that
$q_y=(s+1)f^sf_yM^{(s+1)}+f^{s+1}M^{(s+1)}_y=f^s$. That is,
$$(s+1)f_yM^{(s+1)}+fM^{(s+1)}_y=1~~~~~~~~~~~~~~~(2.36)$$
Let $f(x,y)=f_l(x)y^l+f_{l-1}(x)y^{l-1}+\cdots+f_1(x)y+f_0(x)$,
$M^{(s+1)}=M^{(s+1)}_m(x)y^m+M^{(s+1)}_{m-1}(x)y^{m-1}+\cdots+M^{(s+1)}_1(x)y+M^{(s+1)}_0(x)$
with $f_l(x)M^{(s+1)}_m(x)\neq 0$ and $l\geq 1$. It follows from
equations $(2.32)$ and $(2.30)$ that $f^s_l(x)=tq_t(x)\in{\bf K}^*$.
That is, $f_l(x)\in{\bf K}^*$. It follows from equation $(2.36)$
that
$(s+1)(lf_ly^{l-1}+(l-1)f_{l-1}(x)y^{l-2}+\cdots+f_1(x))(M^{(s+1)}_m(x)y^m+M^{(s+1)}_{m-1}(x)y^{m-1}+\cdots+M^{(s+1)}_1(x)y+M^{(s+1)}_0(x))
+(f_l(x)y^l+f_{l-1}(x)y^{l-1}+\cdots+f_1(x)y+f_0(x))(mM^{(s+1)}_m(x)y^{m-1}+(m-1)M^{(s+1)}_{m-1}(x)y^{m-2}+\cdots+M^{(s+1)}_1(x))=1$~~~~~~~~~~~~~~~~~~~~~~~~~~(2.37)\\
We always view that the polynomials are in ${\bf K}[x][y]$ with
coefficients in ${\bf K}[x]$ when comparing the coefficients of
$y^j$.

Suppose $m\geq 1$. Comparing the coefficients of $y^{l+m-1}$ of
equation $(2.37)$, we have the following equation:
$$(s+1)lf_lM^{(s+1)}_m(x)+mf_lM^{(s+1)}_m(x)=0$$
Since $f_l(x)M^{(s+1)}_m(x)\neq 0$, so we have $(s+1)l+m=0$. This is
a contradiction because $sl\geq 1$ and $m\geq 1$. Therefore, we have
$m=0$. That is, $M^{(s+1)}_y=0$. Then equation $(2.36)$ has the
following form:
$$(s+1)f_yM^{(s+1)}=1$$
Thus, we have $f_y\in{\bf K}^*$. That is, $f=f_1y+f_0(x)$ with
$f_1\in{\bf K}^*$. So $f$ is a polynomial of $y+a(x)$ with
$a(x)=f_1^{-1}f_0(x)$. That is, $q,~u,~h\in{\bf K}[y+a(x)]$. Then
the conclusion follows from the proof of Theorem 2.8 in \cite{13}.\\

Case II If $f(x,y)|q_x$, then
$$q_x=f(x,y)N^{(1)}(x,y)~~~~~~~~~~~~~~(2.38)$$
for some $N^{(1)}(x,y)\in{\bf K}[x,y]$. Thus, it follows from
equations $(2.38)$ and $(2.32)$ that $q_{xy}=f_yN^{(1)}+fN^{(1)}_y$,
$q_{yx}=sf^{s-1}f_x$. Thus, we have the following equation:
$$f_yN^{(1)}=f(sf^{s-2}f_x-N^{(1)}_y).$$
Since $f$ is irreducible and $f_y\neq 0$, so we have $f|N^{(1)}$. That is,
$N^{(1)}=fN^{(2)}$ for some $N^{(2)}\in{\bf K}[x,y]$. It follows
from equation $(2.38)$ that $q_x=f^2N^{(2)}(x,y)$. We can do step by
step until $q_x=f^sN^{(s)}(x,y)$. Then
$q_{xy}=sf^{s-1}f_yN^{(s)}+f^sN^{(s)}_y=q_{yx}=sf^{s-1}f_x$. That
is,
$$s(f_x-f_yN^{(s)})=fN^{(s)}_y~~~~~~~~~~~~~~~(2.39)$$
Since $q_y=f^s$, $q_x=f^sN^{(s)}(x,y)$ and $\deg_yq_x\leq
t-1=\deg_yq_y$, so we have $\deg_y(f^sN^{(s)}(x,y))\leq\deg_y(f^s)$.
Thus, we have $\deg_yN^{(s)}(x,y)=0$. That is, $N^{(s)}_y=0$. Then
equation $(2.39)$ has the following form:
$$f_x=f_yN^{(s)}(x)~~~~~~~~~~~~~~~(2.40)$$
where $N^{(s)}(x):=N^{(s)}(x,0)=N^{(s)}(x,y)$. Let $\bar{x}=x$,
$\bar{y}=y+\int N^{(s)}(x)dx$. Then it follows from equation
$(2.40)$ that $f_{\bar{x}}=0$. That is, $f\in{\bf K}[y+a(x)]$, where
$a(x)=\int N^{(s)}(x)dx\in{\bf K}[x]$. Thus, $q,~u,~h$ are
polynomials of $y+a(x)$. Then the conclusion follows from the proof
of Theorem 2.8 in \cite{13}.
\end{proof}

\begin{cor}
Let $H=(u(x,y),v(x,y,z),h(x,y))$ be a polynomial map over ${\bf
K}[x,y,z]$. Assume that $H(0)=0$ and the components of $H$ are
linearly independent over ${\bf K}$. If $JH$ is nilpotent and the
variety $V(u_y)$ or the variety $V(h_y)$ is irreducible, then $H$
has the form of Theorem 2.2.
\end{cor}
\begin{proof}
Let $v=v_dz^d+\cdots+v_1z+v_0$. Then it follows from Lemma 3.2 in
\cite{14} that $d=1$ and $v_1\in {\bf K}^*$. Since $JH$ is
nilpotent, so we have the following equations:
\begin{equation}
\nonumber
  \left\{ \begin{aligned}
  u_x+v_{0y} = 0~~~~~~~~~~~~~~~~~~~~~~~~(2.4) \\
  u_xv_{0y}-v_{0x}u_y-v_1h_y=0~~~~~~~~(2.5) \\
  v_1(u_xh_y-u_yh_x) = 0~~~~~~~~~~~~~~(2.6)
                          \end{aligned} \right.
  \end{equation}
It follows from equation $(2.6)$ and Lemma 3.1 in \cite{14} that
there exists $q\in{\bf K}[x,y]$ such that
$$u,~h\in{\bf K}[q]~~~~~~~~~~~~~~~~~~~~~~~~~~~~~(2.8)$$
Since $V(u_y)$ or $V(h_y)$ is irreducible and $V(q_y)\subset
V(u_y)\cap V(h_y)$ is a subvariety, so $V(q_y)$ is irreducible. Then
the conclusion follows from the proof of Theorem 2.5.
\end{proof}

\section{Polynomial maps of the form $H=(u(x,y,z),v(x,\allowbreak y,u), h(x,y))$}

In this section, we classify polynomial maps of the form
$H=(u(x,y,z),v(x,y,u),\allowbreak  h(x,y))$ in the case that $JH$ is
nilpotent and $(\deg v(x,y,0),\deg h)\leq 3$. Combining Theorem 2.2
with Corollary 3.7 in \cite{14}, we have the following Proposition.

\begin{prop}
Let $H=(u(x,y,z),v(x,y),h(x,y))$ be a polynomial map over ${\bf
K}[x,y,z]$. Assume that $H(0)=0$ and the components of $H$ are
linearly independent over ${\bf K}$. If $JH$ is nilpotent and
$(\deg_xv, \deg_xh)\leq 3$ or at least one of $\deg_xv$, $\deg_xh$
is a prime, then $u=u_1z-a^{-1}b'(y)g(ax+b(y))-u_1l_2y$,
$v=g(ax+b(y))$, $h=c_0v^2+l_2v$, where
$b(y)=u_1c_0ay^2+l_1y+\tilde{l}_2$, $u_1, c_0, a\in{\bf K}^*$,
$l_1,l_2, \tilde{l}_2\in{\bf K}$, $g(t)\in{\bf K}[t]$ and $g(0)=0$,
$\deg_tg(t)\geq 1$.
\end{prop}

In the following theorem, we denote $v_x=\frac{\partial v}{\partial
x}(x,y,u)$, $v_y=\frac{\partial v}{\partial y}(x,y,u)$ and
$v_u=\frac{\partial v}{\partial u}(x,y,u)$.

\begin{thm}
Let $H=(u(x,y,z),v(x,y,u), h(x,y))$ be a polynomial map over ${\bf
K}[x,y,z]$. Assume that $H(0)=0$ and $\deg_uv\geq 2$. If $JH$ is
nilpotent, then $u,v,h$ are linearly dependent.
\end{thm}
\begin{proof}
If $u_z=0$, then it follows from Proposition 2.1 in \cite{13} that
$u,v,h$ are linearly dependent. Thus, we can assume that $u_z\neq
0$. Since $JH$ is nilpotent, we have the following equations:
\begin{equation}
\nonumber
  \left\{ \begin{aligned}
  u_x+v_y = -v_uu_y~~~~~~~~~~~~~~~~~~~~~~~~(3.1) \\
  u_xv_y-v_xu_y-h_xu_z=h_yv_uu_z~~~~~~~~(3.2) \\
  u_z(v_xh_y-v_yh_x) = 0~~~~~~~~~~~~~~~~~~~(3.3)
                          \end{aligned} \right.
  \end{equation}
Let $u=u_dz^d+u_{d-1}z^{d-1}+\cdots+u_1z+u_0$,
$v(x,y,u)=v_m(x,y)u^m+v_{m-1}(x,y)u^{m-1}+\cdots+v_1(x,y)u+v_0(x,y)$
with $u_dv_m\neq 0$ and $d\geq 1$, $m\geq 2$, $u_i, v_j\in{\bf
K}[x,y]$ for $0\leq i\leq d$, $0\leq j\leq m$. It follows from
equation $(3.1)$ that
$u_{dx}z^d+u_{(d-1)x}z^{d-1}+\cdots+u_{1x}z+u_{0x}+v_{my}u^m+v_{(m-1)y}u^{m-1}+\cdots+v_{1y}u+v_{0y}=
-(mv_mu^{m-1}+\cdots+2v_2u+v_1)(u_{dy}z^d+u_{(d-1)y}z^{d-1}+\cdots+u_{1y}z+u_{0y})$~~~~~~~~~~~~~~(3.4)\\
We always view that the polynomials are in ${\bf K}[x,y][z]$ with
coefficients in ${\bf K}[x,y]$ in the following arguments.\\

Since $m\geq 2$, comparing the coefficients of $z^{md}$ of equation
$(3.4)$, we have $v_{my}u_d+mv_mu_{dy}=0$. That is,
$$\frac{v_{my}}{v_m}=-m\frac{u_{dy}}{u_d}$$
Suppose $v_{my}\neq 0$. Then $u_{dy}\neq 0$. Thus, we have
$v_mu^m_d=e^{c(x)}$ by integrating the two sides of the above
equation with respect to $y$, where $c(x)$ is a function of $x$.
Since $v_m~,u_d\in{\bf K}[x,y]$ and $e^{c(x)}$ is a function of $x$.
Thus, we have $v_m~,u_d\in{\bf K}[x]$. This is a contradiction!
Therefore, we have $v_{my}=0=u_{dy}$. Then we have $u_{iy}=0$ by
comparing the coefficients of $z^{d(m-1)+i}$ of equation $(3.4)$ for
$i=d-1,~d-2,\ldots,1$. Then equation $(3.4)$ has the following
form:\\
$u_{dx}z^d+u_{(d-1)x}z^{d-1}+\cdots+u_{1x}z+u_{0x}+v_{(m-1)y}u^{m-1}+\cdots+v_{1y}u+v_{0y}=
-(mv_mu^{m-1}+\cdots+2v_2u+v_1)u_{0y}$~~~~~~~~~~~~~~~~~~~~~~~~~~~~~~~~~~~~~~~~(3.5)\\
Comparing the coefficients of $z^{jd}$ for $j=m-1,~m-2,\ldots,2$ of
equation $(3.5)$, we have the following equations:
$$v_{jy}+(j+1)v_{j+1}u_{0y}=0~~~~~~~~~~~~~~~~~~~(3.6)$$
Then we have
$$(v_{1y}+2v_2u_{0y})u_k=-u_{kx}~~~~~~~~~~~~~~~~~(3.7)$$
by comparing the coefficients of $z^k$ of equation $(3.5)$ for
$k=d,~d-1,\ldots,1$. If $v_{1y}+2v_2u_{0y}=0$, then $u_{kx}=0$ for
$1\leq k\leq d$. If $v_{1y}+2v_2u_{0y}\neq 0$, then we have
$u_{kx}=0$ by comparing the degree of $x$ of two sides of equation
$(3.7)$ for $1\leq k\leq d$. Since $u_d\neq 0$, we have
$v_{1y}+2v_2u_{0y}=0$. This is a contradiction! Thus, we have
$v_{1y}+2v_2u_{0y}=0$ and $u_{dx}=\cdots=u_{1x}=0$. Therefore, we
have $u_d,\ldots,u_1\in{\bf K}^*$. Then equation $(3.5)$ has the
following form:
$$u_{0x}+v_{0y}+v_1u_{0y}=0~~~~~~~~~~~~~~~~~~~~~~~(3.8)$$

Since $u_z\neq 0$, so it follows from equation (3.3) that
$h_yv_x-h_xv_y=0$. That is,
$$h_y(v_{mx}u^m+\cdots+v_{1x}u+v_{0x})-h_x(v_{(m-1)y}u^{m-1}+\cdots+v_{1y}u+v_{0y})=0~~~~~(3.9)$$
Then we have $h_yv_{mx}=0$ by comparing the coefficients of $z^{md}$
of equation $(3.9)$. Thus, we have $h_y=0$ or $v_{mx}=0$.\\

(1) If $h_y=0$, then it follows from equation $(3.9)$ that $h_x=0$
or $v_{(m-1)y}=\cdots=v_{1y}=v_{0y}=0$.

(i) If $h_x=0$, then $h=0$ because $h(0)=0$. Thus, $u,v,h$ are
linearly dependent.

(ii) If $v_{(m-1)y}=\cdots=v_{0y}=0$, then $v_y=0$. It follows from
equation $(3.2)$ that $u_yv_x+h_xu_z=0$. That is,
$u_{0y}(v_{mx}u^m+v_{(m-1)x}u^{m-1}+\cdots+v_{1x}u+v_{0x})+h_x(du_dz^{d-1}+(d-1)u_{d-1}z^{d-2}+\cdots+u_1)=0$
~~~~~~~~~~~~~~~~~~~~~(3.10)

If $u_{0y}=0$, then it follows from equation $(3.10)$ that $h_x=0$.
Since $H(0)=0$, we have $h=0$. Thus, $u,v,h$ are linearly dependent.

If $u_{0y}\neq 0$, then we have $v_{mx}=\cdots=v_{1x}=0$ by
comparing the coefficients of
$z^{md},z^{(m-1)d},\ldots,z^d$ of equation $(3.10)$ respectively.\\
If $d\geq 2$, then we have $h_x=0$ by comparing the coefficient of
$z^{d-1}$ of equation $(3.10)$.  Since $H(0)=0$, we have $h=0$.
Thus, $u,v,h$ are linearly dependent. \\
If $d=1$, then it follows from equation $(3.6)$ ($j=m-1$) that $u_{0y}=0$. This is a contradiction! \\

(2) If $v_{mx}=0$, then $v_m\in {\bf K}^*$. It follows from equation
$(3.2)$ that\\
$u_{0x}(v_{(m-1)y}u^{m-1}+\cdots+v_{1y}u+v_{0y})-u_{0y}(v_{(m-1)x}u^{m-1}+\cdots+v_{1x}u+v_{0x})-h_x(du_dz^{d-1}+(d-1)u_{d-1}z^{d-2}+\cdots+u_1)
-h_y(mv_mu^{m-1}+\cdots+2v_2u+v_1)(du_dz^{d-1}+(d-1)u_{d-1}z^{d-2}+\cdots+u_1)=0$~~~~~~~~~~~~~~~~~~~~~~~(3.11)

If $d\geq 2$, then we have $h_y=0$ by comparing the coefficient of
$z^{d(m-1)+d-1}$ of equation $(3.11)$. Thus, we have
$u_{0x}v_{ly}-u_{0y}v_{lx}=0$ by comparing the coefficients of
$z^{ld}$ of equation $(3.11)$ for $l=m-1,m-2,\ldots,1$. Comparing
the coefficients of $z^{d-1}$ of equation $(3.11)$, we have $h_x=0$.
Thus, we have $h=0$ because $h(0)=0$. So $u,v,h$ are linearly
dependent.

If $d=1$, then equation $(3.11)$ has the following form:\\
$u_{0x}(v_{(m-1)y}u^{m-1}+\cdots+v_{1y}u+v_{0y})-u_{0y}(v_{(m-1)x}u^{m-1}+\cdots+v_{1x}u+v_{0x})-h_xu_1
-h_yu_1(mv_mu^{m-1}+\cdots+2v_2u+v_1)=0$~~~~~~~~~~~~~~~~~~~~~~~~~~~~~~~~~~(3.12)\\
Thus, we have
$$u_{0x}v_{ly}-u_{0y}v_{lx}-h_yu_1(l+1)v_{l+1}=0~~~~~~~~~~~~~~~~~~~~~(3.13)$$
by comparing the coefficients of $z^{ld}$ of equation $(3.12)$ for
$l=m-1,m-2,\ldots,1$. Then equation $(3.12)$ has the following form:
$$u_{0x}v_{0y}-u_{0y}v_{0x}-u_1(h_x+v_1h_y)=0~~~~~~~~~~~~~~~~~~~~~~~~(3.14)$$
Comparing the coefficients of $z^{ld}$ of equation $(3.9)$ for
$l=m-1,m-2,\ldots,1$, we have
$$h_yv_{lx}-h_xv_{ly}=0~~~~~~~~~~~~~~~~~~~~~~~~~~~~~(3.15)$$
for $l=m-1,\ldots,1$. Then equation $(3.9)$ has the following form:
$$h_yv_{0x}-h_xv_{0y}=0~~~~~~~~~~~~~~~~~~~~~~~~~~~~~(3.16)$$
It follows from equation $(3.6)$ that
$$v_{m-1}=-mv_mu_0+c_1(x)~~~~~~~~~~~~~~~~~~~~~~~~~~~(3.17)$$
for some $c_1(x)\in{\bf K}[x]$. Substituting equation (3.17) to
equation (3.13) for $l=m-1$, we have
$$u_{0y}c'_1(x)+h_ymv_mu_1=0$$
Thus, we have $u_0c'_1(x)+mv_mu_1\cdot h+e(x)=0$ by integrating the
two sides of the above equation with respect to $y$, where $e(x)\in
{\bf K}[x]$. That is,
$$h=-\frac{1}{mv_mu_1}c'_1(x)u_0-\frac{1}{mv_mu_1}e(x)~~~~~~~~~~~~~~~(3.18)$$
Substituting equations $(3.17)$ and $(3.18)$ to equation $(3.15)$
for $l=m-1$, we have
$$c''_1(x)u_0u_{0y}+\frac{1}{mv_m}(c'_1(x))^2u_{0y}+e'(x)u_{0y}=0$$
If $u_{0y}=0$, then it follows from equation $(3.18)$ that $h_y=0$.
Then it reduces to (1).\\
If $u_{0y}\neq 0$, then we have
$$c''_1(x)u_0+\frac{1}{mv_m}(c'_1(x))^2+e'(x)=0~~~~~~~~~~~~(3.19)$$
Thus, we have $c_1''(x)=0$ and
$$e'(x)=-\frac{1}{mv_m}(c_1'(x))^2\in{\bf K}~~~~~~~~~~~~~~~~~~~(3.20)$$
by comparing the degree of $y$ of the equation (3.19). Then we have
$$c_1'(x)(u_{0x}v_{0y}-u_{0y}v_{0x})+e'(x)v_{0y}=0~~~~~~~~~(3.21)$$
by substituting equation $(3.18)$ to equation $(3.16)$. Substituting
equation (3.20) to equation (3.21), we have
$$c_1'(x)(u_{0x}v_{0y}-u_{0y}v_{0x})-\frac{1}{mv_m}(c_1'(x))^2v_{0y}=0~~~~~~~~~(3.22)$$

If $c_1'(x)=0$, then it follows from equation (3.20) that
$e(x)\in{\bf K}$. Thus, it follows from equation (3.18) that
$h=-\frac{e(x)}{mv_mu_1}\in {\bf K}$. Since $h(0,0)=0$, we have
$h=0$. Therefore, $u,v,h$ are linearly dependent.

If $c_1'(x)\neq 0$, then equation (3.22) has the following form:
$$u_{0x}v_{0y}-u_{0y}v_{0x}=\frac{1}{mv_m}c_1'(x)v_{0y}~~~~~~~~~~~~~~~~~(3.23)$$
Substituting equations (3.18) and (3.23) to equation (3.14), we have
the following equation:
$$c_1'(x)(v_{0y}+u_{0x}+v_1u_{0y})+e'(x)=0.$$
Substituting equation (3.8) to the above equation, we have
$e'(x)=0$. It follows from equation (3.20) that $c_1'(x)=0$. This is
a contradiction!
\end{proof}

\begin{thm}
Let $H=(u(x,y,z),v(x,y,u), h(x,y))$ be a polynomial map over ${\bf
K}[x,y,z]$. Assume that $H(0)=0$ and the components of $H$ are
linearly independent over ${\bf K}$. If $JH$ is nilpotent and $(\deg
v(x,y,0), \deg h)\leq 3$, then there exists $T\in
\operatorname{GL}_3({\bf K})$ such that $T^{-1}HT$ has the form of
Proposition 3.1.
\end{thm}
\begin{proof}
It follows from Theorem 3.2 that $\deg_uv\leq 1$. If $u_z=0$ or
$v_u=0$, then it reduces to Proposition 2.1 in \cite{13} and
Proposition 3.1 respectively. Thus, we can assume that $u_z\neq 0$
and $v_u\neq 0$. Since $JH$ is nilpotent, we have the following
equations:
\begin{equation}
\nonumber
  \left\{ \begin{aligned}
  u_x+v_y = -v_uu_y~~~~~~~~~~~~~~~~~~~~~~~~(3.1) \\
  u_xv_y-v_xu_y-h_xu_z=h_yv_uu_z~~~~~~~~(3.2) \\
  u_z(v_xh_y-v_yh_x) = 0~~~~~~~~~~~~~~~~~~~(3.3)\\
                          \end{aligned} \right.
  \end{equation}
Let $u=u_dz^d+u_{d-1}z^{d-1}+\cdots+u_1z+u_0$ and $v=v_1u+v_0$ with
$u_dv_1\neq 0$.

It follows from equation (3.1) that
$u_{dx}z^d+u_{(d-1)x}z^{d-1}+\cdots+u_{1x}z+u_{0x}+v_{1y}u+v_{0y}+v_1(u_{dy}z^d+u_{(d-1)y}z^{d-1}+\cdots+u_{1y}z+u_{0y})=0$
~~~~~~~~~~~~~~(3.24)\\
Then we have
$$u_{dx}+v_{1y}u_d+v_1u_{dy}=0~~~~~~~~~~~~~~~~~~~~(3.25)$$
and
$$u_{(d-1)x}+v_{1y}u_{d-1}+v_1u_{(d-1)y}=0~~~~~~~~~~(3.26)$$
by comparing the coefficients of $z^d$ and $z^{d-1}$ of equation
(3.24) respectively.

It follows from equation (3.2) that\\
$(u_{dx}z^d+u_{(d-1)x}z^{d-1}+\cdots+u_{1x}z+u_{0x})(v_{1y}u+v_{0y})-(u_{dy}z^d+u_{(d-1)y}z^{d-1}+\cdots+u_{1y}z+u_{0y})(v_{1x}u+v_{0x})
-h_x(du_dz^{d-1}+(d-1)u_{d-1}z^{d-2}+\cdots+u_1)-h_yv_1(du_dz^{d-1}+(d-1)u_{d-1}z^{d-2}+\cdots+u_1)=0$~~~~~~~~~~~~~~~~~~~~~~~~~~~~~~~~~~~(3.27)\\
Comparing the coefficients of $z^k$ of equation (3.27) for
$k=2d,2d-1,\ldots,d+1$, we have
$$u_{(k-d)x}v_{1y}-u_{(k-d)y}v_{1x}=0~~~~~~~~~~~~~~~~~~~~~~(3.28)$$
Then we have
$$(u_{0x}v_{1y}-u_{0y}v_{1x})u_d+u_{dx}v_{0y}-u_{dy}v_{0x}=0~~~~~~~~~(3.29)$$
and
$$(u_{0x}v_{1y}-u_{0y}v_{1x})u_{d-1}+u_{(d-1)x}v_{0y}-u_{(d-1)y}v_{0x}-du_d(h_x+v_1h_y)=0~~~~~~(3.30)$$
by comparing the coefficients of $z^d$ and $z^{d-1}$ of equation
(3.27) respectively.

It follows from equation (3.3) that
$$h_y(v_{1x}u+v_{0x})-h_x(v_{1y}u+v_{0y})=0~~~~~~~~~~~~~~~~(3.31)$$
Then we have
$$h_yv_{1x}-h_xv_{1y}=0~~~~~~~~~~~~~~~~~~~~~~~~~~~~~~(3.32)$$
and
$$h_yv_{0x}-h_xv_{0y}=0~~~~~~~~~~~~~~~~~~~~~~~~~~~~~~(3.33)$$
by comparing the coefficients of $z^d$ and $z^0$ of equation (3.31)
orderly.\\

If $h_y\neq 0$, then let $\lambda(x,y)=\frac{h_x}{h_y}$. Thus, we
have $h_x=\lambda(x,y)h_y$. It follows from equations (3.32) and
(3.33) that
$$v_{1x}=\lambda(x,y)v_{1y}~~~~~~~~~~~~~~~~~(3.34)$$
and
$$v_{0x}=\lambda(x,y)v_{0y}~~~~~~~~~~~~~~~~~(3.35)$$
It follows from equations (3.28) and (3.34) that
$$(u_{(k-d)x}-\lambda(x,y)u_{(k-d)y})v_{1y}=0$$

If $v_{1y}=0$, then it follows from equation (3.32) that $v_{1x}=0$
because $h_y\neq 0$. That is, $v_1\in{\bf K}^*$.

If $v_{1y}\neq 0$, then we have
$$u_{(k-d)x}=\lambda(x,y)u_{(k-d)y}~~~~~~~~~~~~~~~(3.36)$$
for $d+1\leq k\leq 2d$. Substituting equations (3.34), (3.35),
(3.36)$(k=2d)$ to equation (3.29), we have the following equation:
$$u_{0x}=\lambda(x,y)u_{0y}~~~~~~~~~~~~~~~~(3.37)$$
Substituting equations (3.34), (3.35), (3.36), (3.37) to equation
(3.30), we have the following equation:
$$h_x+v_1h_y=0.$$
That is,
$$v_1=-\frac{h_x}{h_y}=-\lambda(x,y).$$
Substituting equation (3.34) to the above equation, we have
$$v_1v_{1y}=-v_{1x}.$$
Then we have $v_1\in{\bf K}^*$ by comparing the degree of $x$ of two
sides of the above equation. This is a contradiction!\\
Consequently, we have $H=(u,v_1u+v_0(x,y),h(x,y))$ with $v_1\in{\bf
K}^*$. Let
\[T=\left(
  \begin{array}{ccc}
    1 & 0 & 0 \\
    v_1 & 1 & 0 \\
    0 & 0 & 1 \\
  \end{array}
\right).\] Then $T^{-1}HT=(u(x,y+v_1x,z),v_0(x,y+v_1x),h(x,y+v_1x))$
and $\deg_xv_0(x,y+v_1x)=\deg v(x,y,0)$, $\deg_xh(x,y+v_1x)=\deg
h(x,y)$. Since $(\deg v(x,y,0), \deg h)\leq 3$, so the conclusion
follows from Proposition 3.1.\\

If $h_y=0$, then it follows from equation (3.32) that $h_x=0$ or
$v_{1y}=0$.

(1) If $h_x=0$, then $h=0$ because $h(0,0)=0$. Thus, $u,v,h$ are
linearly dependent. This is a contradiction!

(2) If $v_{1y}=0$, then it follows from equation $(3.28)$ that
$v_{1x}=0$ or $u_{dy}=u_{(d-1)y}=\cdots=u_{1y}=0$.

Case I If $v_{1x}=0$, then $v_1\in{\bf K}^*$. Thus, the conclusion
follows from the former arguments.

Case II If $u_{dy}=u_{(d-1)y}=\cdots=u_{1y}=0$, then it follows from
equations (3.25) and (3.26) that $u_{dx}=0=u_{(d-1)x}$. Thus, it
follows from equation (3.29) that $v_{1x}=0$ or $u_{0y}=0$. If
$v_{1x}=0$, then we have $v_1\in{\bf K}^*$, it reduces to Case I. If
$u_{0y}=0$, then it follows from equation (3.30) that $h_x=0$. Thus,
we have $h=0$ because $h(0,0)=0$. Therefore, $u,v,h$ are linearly
dependent. This is a contradiction!
\end{proof}

\begin{cor}
Let $H=(u(x,y,v),v(x,y,z),h(x,y))$ be a polynomial map with
$H(0)=0$. Assume that the components of $H$ are linearly independent
over ${\bf K}$. If $JH$ is nilpotent and $(\deg u(x,y,0),\deg
h(x,y))\leq 3$, then there exists $T\in \operatorname{GL}_3({\bf
K})$ such that $T^{-1}HT$ has the form of Proposition 3.1.
\end{cor}
\begin{proof}
Let
\[T_1=\left(
  \begin{array}{ccc}
    0 & 1 & 0 \\
    1 & 0 & 0 \\
    0 & 0 & 1 \\
  \end{array}
\right).\] Then $T_1^{-1}HT_1=(v(y,x,z),u(y,x,v(y,x,z)),h(y,x))$.
Since $JH$ is nilpotent, we have that
$J(T_1^{-1}HT_1)=T_1^{-1}JHT_1$ is nilpotent. Since $\deg
u(y,x,0)=\deg u(x,y,0)$, $\deg h(y,x))=\deg h(x,y)$ and $(\deg
u(x,y,0),\deg h(x,y))\leq 3$, so it follows from Theorem 3.3 that
there exists $T_2\in \operatorname{GL}_3({\bf K})$ such that
$T_2^{-1}(T_1^{-1}HT_1)T_2$ is of the form of Proposition 3.1. Let
$T=T_1T_2\in \operatorname{GL}_3({\bf K})$. Then the conclusion
follows.
\end{proof}

\section{Polynomial maps of the form $H=(u(x,y,z),v(x,\allowbreak y,z), h(x,y))$}

In the section, we classify polynomial maps of the form
$H=(u(x,y,z),v(x,y,z),\allowbreak h(x,y))$ in the case that $JH$ is
nilpotent and $\deg_zv(x,y,z)\leq 1$. Firstly, we prove that $u,v,h$
are linearly dependent in the case that $JH$ is nilpotent and
$\deg_zv(x,y,z)=1$ and $\deg_zu\geq 2$. In the proof of the
following theorem, we divide two cases according to the degree of
$z$ in $u$. Case I: $\deg_zu=2$; Case II: $\deg_zu\geq 3$. In the
first case, we divide two parts. We have three subcases in the
second case.

\begin{thm}
Let $H=(u(x,y,z),v(x,y,z),h(x,y))$ be a polynomial map with
$\deg_zv(x,y,z)=1$. Assume that $H(0)=0$ and $\deg_zu\geq 2$. If
$JH$ is nilpotent, then $u,v,h$ are linearly dependent.
\end{thm}
\begin{proof}
Let $u=u_dz^d+\cdots+u_1z+u_0$, $v=v_1z+v_0$ with $u_dv_1\neq 0$ and
$d\geq 2$. Since $JH$ is nilpotent, we have the following equations:
\begin{equation}
\nonumber
  \left\{ \begin{aligned}
  u_x+v_y = 0~~~~~~~~~~~~~~~~~~~~~~~~~~~~~~~~~~~~~~~~~~~~~~~(4.1) \\
  u_xv_y-v_xu_y-h_xu_z-h_yv_z=0~~~~~~~~~~~~~~~~~~~~~~(4.2) \\
  v_xh_yu_z-h_xv_yu_z+h_xu_yv_z-u_xh_yv_z = 0~~~~~~~~~~~(4.3)\\
                          \end{aligned} \right.
  \end{equation}
It follows from equation (4.1) that
$$u_{dx}z^d+u_{(d-1)x}z^{d-1}+\cdots+u_{1x}z+u_{0x}+v_{1y}z+v_{0y}=0~~~~~~~~~~~(4.4)$$
We always view that the polynomials are in ${\bf K}[x,y][z]$ with
coefficients in ${\bf K}[x,y][z]$ in the following arguments. Thus,
we have
$$u_{dx}=\cdots=u_{2x}=0~~~~~~~~~~~~~~~(4.5)$$
and
$$u_{ix}+v_{iy}=0~~~~~~~~~~~~~~~~~~~~~~~(4.6)$$
for $i=1,~0$ by comparing the coefficients of $z^j$ of equation
$(4.4)$ for $j=d,\ldots,1,0$.

It follows from equations (4.2) and (4.5) that\\
$(u_{1x}z+u_{0x})(v_{1y}z+v_{0y})-(v_{1x}z+v_{0x})(u_{dy}z^d+u_{(d-1)y}z^{d-1}+\cdots+u_{1y}z+u_{0y})-h_x(du_dz^{d-1}+\cdots+u_1)-h_yv_1=0$
~~~~~~~~~~~~~~~~~~~~~~~~~(4.7)\\

Case I If $d=2$, then comparing the coefficients of $z^3$ and $z^2$
of equation (4.7), we have $v_{1x}u_{2y}=0$ and
$$u_{1x}v_{1y}-v_{1x}u_{1y}-v_{0x}u_{2y}=0~~~~~~~~~~~~(4.8)$$
Thus, we have $u_{2y}=0$ or $v_{1x}=0$.\\

(1) If $u_{2y}=0$, then it follows from equations (4.5) and (4.8)
that $u_2\in{\bf K}^*$ and
$$u_{1x}v_{1y}-v_{1x}u_{1y}=0~~~~~~~~~~~~~~~~~~~~~(4.9)$$
It follows from equation (4.6) that $u_{1x}=-v_{1y}$, so there
exists $P\in{\bf K}[x,y]$ such that $u_1=P_y$, $v_1=-P_x$. It
follows from equation (4.9) that $P_{xy}^2-P_{xx}P_{yy}=0$. Then it
follows from Lemma 2.1 in \cite{18} that
$$u_1=P_y=bf(ax+by)+c_2~~~~~~~~~~~~~~~~~~~~~(4.10)$$
and
$$v_1=-P_x=-af(ax+by)+c_1~~~~~~~~~~~~~~~~~~(4.11)$$
for some $f(t)\in {\bf K}[t]$ and $f(0)=0$, $a,b\in{\bf K}^*$,
$c_1,c_2\in{\bf K}$. Then equation (4.7) has the following form:
$$(u_{1x}z+u_{0x})(v_{1y}z+v_{0y})-(v_{1x}z+v_{0x})(u_{1y}z+u_{0y})-h_x(2u_2z+u_1)-h_yv_1=0~~~(4.12)$$
Thus, we have the following equations
$$u_{1x}v_{0y}+u_{0x}v_{1y}-v_{1x}u_{0y}-v_{0x}u_{1y}-2u_2h_x=0~~~~~~~~~~(4.13)$$
and
$$u_{0x}v_{0y}-v_{0x}u_{0y}-h_xu_1-h_yv_1=0~~~~~~~~~~~~~~~~~~~~~~~~~(4.14)$$
by comparing the coefficients of $z$ and $z^0$ of equation (4.12).
It follows from equations (4.3) and (4.5) that
$$[(v_{1x}z+v_{0x})h_y-h_x(v_{1y}z+v_{0y})](2u_2z+u_1)+[h_x(u_{1y}z+u_{0y})-h_y(u_{1x}z+u_{0x})]v_1=0~~~(4.15)$$
Comparing the coefficients of $z^2,z,z^0$ of equation $(4.15)$, we
have the following equations:
$$v_{1x}h_y-h_xv_{1y}=0~~~~~~~~~~~~~~~~~~~~~~~~~~~~~~~~~~~~~~~~~~(4.16)$$
and
$$(v_{0x}h_y-h_xv_{0y})2u_2+(v_{1x}h_y-h_xv_{1y})u_1+(h_xu_{1y}-h_yu_{1x})v_1=0~~~~~~(4.17)$$
and
$$(v_{0x}h_y-h_xv_{0y})u_1+(h_xu_{0y}-h_yu_{0x})v_1=0~~~~~~~~~~~~~~~~~~~~~~~~~~(4.18)$$
It follows from equations (4.11) and (4.16) that $af'\cdot
(bh_x-ah_y)=0$. Thus, we have $a=0$ or $f'=0$ or $bh_x=ah_y$.

(i) If $a=0$, then
$$u_1=bf(by)+c_2,~~~v_1=c_1\in {\bf K}^*~~~~~~~~~~~~~~~~~(4.19)$$
It follows from equation (4.13) that $2u_2h_x=-v_{0x}u_{1y}$. That
is,
$$2u_2h_x=-b^2f'(by)v_{0x}.$$
Integrating with respect to $x$ of two sides of the above equation,
we have
$$h=-\frac{b^2}{2u_2}f'(by)v_0+\frac{c(y)}{2u_2}~~~~~~~~~~~~~~~~~~~~~~~~~~(4.20)$$
for some $c(y)\in {\bf K}[y]$. Substituting equations (4.19) and
(4.20) to equation (4.17), we have the following equation:
$$v_{0x}[-b^3f''(by)v_0+c'(y)-\frac{b^4}{2u_2}v_1\cdot (f'(by))^2]=0.$$
Thus, we have $v_{0x}=0$ or $f''(by)=0$ and
$c'(y)=\frac{b^4}{2u_2}v_1\cdot (f'(by))^2$.

If $v_{0x}=0$, then it follows from equations (4.13) and (4.19) that
$h_x=0$. It follows from equation (4.18) that $v_1h_yu_{0x}=0$.
Thus, we have $h_y=0$ or $u_{0x}=0$. If $u_{0x}=0$, then it follows
from equation (4.14) that $h_y=0$. If $h_y=0$, then $h=0$ because
$h(0,0)=0$. Thus, $u,v,h$ are linearly dependent.

If $f''(by)=0$, then $f'(by)\in {\bf K}$ and
$c'(y)=\frac{b^4}{2u_2}v_1(f'(by))^2\in {\bf K}$. Let $l:=b^2f'(by)$
and $c:=\frac{l^2}{(2u_2)^2}v_1$. Then it follows from equation
(4.19) that $v_1\in {\bf K}^*$ and $u_1=ly+c_2$ for some $c_2\in{\bf
K}$. Since $h(0,0)=0$, so it follows from equation (4.20) that
$$h=-\frac{l}{2u_2}v_0+c\cdot y~~~~~~~~~~~~~~~~~~~~~~~~~~~~~~~~(4.21)$$
It follows from equations (4.14) and (4.21) that
$$v_{0x}u_{0y}-u_{0x}v_{0y}=\frac{l}{2u_2}u_1v_{0x}+\frac{l}{2u_2}v_1v_{0y}-cv_1~~~~~~~~(4.22)$$
It follows from equations (4.18) and (4.21) that
$$c\cdot u_1v_{0x}+v_1[-\frac{l}{2u_2}(v_{0x}u_{0y}-u_{0x}v_{0y})-c\cdot u_{0x}]=0~~~~~~~~~~~(4.23)$$
Substituting equation (4.22) to equation (4.23), we have the
following equation:
$$c\cdot u_1v_{0x}+v_1[-\frac{l^2}{(2u_2)^2}u_1v_{0x}+\frac{lc}{2u_2}v_1-\frac{l^2}{(2u_2)^2}v_1v_{0y}-c\cdot u_{0x}]=0.$$
Since $c=\frac{l^2}{(2u_2)^2}v_1$, so the above equation has the
following form:
$$\frac{lc}{2u_2}v_1-c(v_{0y}+u_{0x})=0.$$
Substituting equation (4.6)$(i=0)$ to the above equation, we have
$\frac{lc}{2u_2}v_1=0$. That is, $lc=0$. Since
$c=\frac{l^2}{(2u_2)^2}v_1$, so we have $c=l=0$. It follows from
equation (4.21) that $h=0$. Thus, $u,v,h$ are linearly dependent.

(ii) If $f'=0$, then $f=0$ because $f(0)=0$. That is, $u_1=c_2$,
$v_1=c_1\in {\bf K}^*$. It follows from equation (4.13) that
$h_x=0$. It follows from equation (4.17) that $v_{0x}h_y=0$. Thus,
we have $h_y=0$ or
$v_{0x}=0$.\\
If $h_y=0$, then $h=0$ because $h(0,0)=0$. Thus, $u,v,h$ are
linearly dependent.\\
If $v_{0x}=0$, then it follows from equation (4.18) that
$u_{0x}h_y=0$. That is, $u_{0x}=0$ or $h_y=0$. If $u_{0x}=0$, then
it follows from equation (4.14) that $h_y=0$. Thus, we have that
$h_y=0$. It reduces to the above case.

(iii) If $bh_x=ah_y$, then we have
$$h_y=\frac{b}{a}h_x~~~~~~~~~~~~~~~~~~~~(4.24)$$
because we can assume that $a\cdot f'\neq 0$. Otherwise, it reduces
to the former two cases. Let $\bar{x}=ax+by$, $\bar{y}=y$. Then we
have that $h_{\bar{y}}=0$. That is, $h\in {\bf K}[ax+by]$. It
follows from equations (4.10), (4.11), (4.17),(4.24) that
$$h_x(\frac{b}{a}v_{0x}-v_{0y})=0~~~~~~~~~~~~~~~(4.25)$$
It follows from equations (4.18) and (4.24) that
$$u_1h_x(\frac{b}{a}v_{0x}-v_{0y})+v_1h_x(u_{0y}-\frac{b}{a}u_{0x})=0~~~~~~~(4.26)$$
It follows from equations (4.25) and (4.26) that $h_x=0$ or
$bv_{0x}=av_{0y}$ and $bu_{0x}=au_{0y}$. \\
If $h_x=0$, then it follows from equation (4.24) that $h_y=0$. Thus,
we have $h=0$ because $h(0,0)=0$. so $u,v,h$ are linearly
dependent.\\
If $bv_{0x}=av_{0y}$ and $bu_{0x}=au_{0y}$, then $v_0,u_0\in {\bf
K}[ax+by]$. Thus, it follows from equations (4.10), (4.11) and
(4.13) that $h_x=0$. It reduces to the former case.\\

(2) If $v_{1x}=0$, then it follows from equation (4.8) that
$$u_{1x}v_{1y}-v_{0x}u_{2y}=0~~~~~~~~~~~~~~~~~(4.27)$$
It follows from equation (4.7) that
$$(u_{1x}z+u_{0x})(v_{1y}z+v_{0y})-v_{0x}(u_{2y}z^2+u_{1y}z+u_{0y})-h_x(2u_2z+u_1)-h_yv_1=0~~~~(4.28)$$
Comparing the coefficients of $z^2,~z,~z^0$ of equation (4.28), we
have the following equations:
\begin{equation}
\nonumber
  \left\{ \begin{aligned}
  u_{1x}v_{1y}-v_{0x}u_{2y} = 0~~~~~~~~~~~~~~~~~~~~~~~~~~~~~~~(4.29) \\
  u_{1x}v_{0y}+u_{0x}v_{1y}-v_{0x}u_{1y}-2u_2h_x = 0~~~~~~~~(4.30) \\
  u_{0x}v_{0y}-v_{0x}u_{0y}-u_1h_x-v_1h_y = 0~~~~~~~~~~~~(4.31)\\
                          \end{aligned} \right.
  \end{equation}
It follows from equations (4.3) and (4.5) that
$$[v_{0x}h_y-(v_{1y}z+v_{0y})h_x](2u_2z+u_1)+[h_x(u_{2y}z^2+u_{1y}z+u_{0y})-h_y(u_{1x}z+u_{0x})]v_1=0$$
Comparing the coefficients of $z^2,~z$ and $z^0$ of the above
equation, we have the following equations:
\begin{equation}
\nonumber
  \left\{ \begin{aligned}
  h_x(v_1u_{2y}-2u_2v_{1y}) = 0~~~~~~~~~~~~~~~~~~~~~~~~~~~~~~~~~~~~~~~~~~~~~(4.32) \\
  -v_{1y}h_xu_1+2u_2(v_{0x}h_y-v_{0y}h_x)+v_1(h_xu_{1y}-h_yu_{1x}) = 0~~~~(4.33) \\
  u_1(v_{0x}h_y-v_{0y}h_x)+v_1(h_xu_{0y}-h_yu_{0x}) = 0~~~~~~~~~~~~~~~~~~~(4.34)\\
                          \end{aligned} \right.
  \end{equation}
It follows from equation (4.32) that $h_x=0$ or
$v_1u_{2y}=2u_2v_{1y}$.

If $h_x=0$, then it follows from equation (4.33) that
$h_y(2u_2v_{0x}-u_{1x}v_1)=0$. Thus, we have that $h_y=0$ or
$2u_2v_{0x}=v_1u_{1x}$. If $h_y=0$, then $h=0$ because $h(0,0)=0$.
Thus, $u,v,h$ are linearly dependent. If $2u_2v_{0x}=v_1u_{1x}$,
then it follows from equation (4.6) that
$$2u_2v_{0x}=-v_1v_{1y}\in {\bf K}[y]~~~~~~~~~~~~~(4.35)$$
Substituting equations (4.35) and (4.6) to equation (4.29), we have
the following equation: $v_{1y}(2u_2v_{1y}-v_1u_{2y})=0$. Thus, we
have $v_{1y}=0$ or $2u_2v_{1y}=v_1u_{2y}$. \\
If $v_{1y}=0$, then it follows from equation (4.29) that
$v_{0x}u_{2y}=0$. Thus, we have $v_{0x}=0$ or $u_{2y}=0$. If
$u_{2y}=0$, then it reduces to Case I(1). If $v_{0x}=0$, then it
follows form equation (4.34) that $h_yu_{0x}=0$. Thus, we have
$h_y=0$ or $u_{0x}=0$. If $u_{0x}=0$, then it follows from equation
(4.31) that $h_y=0$. Therefore, we have $h=0$ because $h(0,0)=0$. So
$u,v,h$ are linearly dependent.

If $2u_2v_{1y}=v_1u_{2y}$, then we have
$$\frac{u_{2y}}{u_2}=2\frac{v_{1y}}{v_1}~~~~~~~~~~~~~~~~~~~~(4.36)$$
Suppose $u_{2y}v_{1y}\neq 0$. Then we have
$$u_2=e^{\bar{c}(x)}v_1^2$$
by integrating of two sides of equation (4.36) with respect to $y$.
where $\bar{c}(x)$ is a function of $x$. Since $u_2,v_1\in {\bf
K}[x,y]$, we have $e^{\bar{c}(x)}\in {\bf K}[x]$. That is,
$u_2=c(x)v_1^2$ for $c(x)\in {\bf K}[x]$ and $c(x)\neq 0$. Then it
follows from equations (4.29) and (4.6) that
$$v_{1y}(2c(x)v_1v_{0x}+v_{1y})=0.$$
That is,
$$2c(x)v_1v_{0x}=-v_{1y}~~~~~~~~~~~~~~~~~~(4.37)$$
If $v_{0x}\neq 0$, then we have that $v_{1y}=0$ by comparing the
degree of $y$ of equation (4.37). Thus, we have $v_{0x}=0$. This is
a contradiction. Therefore, we have $v_{1y}=v_{0x}=0$. It follows
from equation (4.36) that $u_{2y}=0$, This is a contradiction! If
$u_{2y}v_{1y}=0$, then it follows from equation (4.36) that
$u_{2y}=v_{1y}=0$. which reduces to case I(1).\\

Case II If $d\geq 3$, then we have
$$u_{dy}v_{1x}=0~~~~~~~~~~~~~~~~~~~~~~(4.38)$$
and
$$v_{1x}u_{(d-1)y}+v_{0x}u_{dy}=0~~~~~~~~~~~~~~~~~~(4.39)$$
by comparing the coefficients of $z^{d+1}$ and $z^d$ of equation
(4.7). It follows from equations (4.38) and (4.39) that
$v_{1x}=v_{0x}=0$ or $v_{1x}=0=u_{dy}$ or $u_{dy}=u_{(d-1)y}=0$.\\

(a) If $v_{1x}=v_{0x}=0$, then equation (4.7) has the following
form:
$$(u_{1x}z+u_{0x})(v_{1y}z+v_{0y})-h_x(du_dz^{d-1}+\cdots+u_1)-h_yv_1=0~~~~~~~~~~~~~~~~~(4.40)$$
If $d>3$, then $h_x=0$ by comparing the coefficient of $z^{d-1}$ of
equation (4.40). Thus, it follows from equation (4.3) that
$h_y(u_{1x}z+u_{0x})=0$. Therefore, we have $h_y=0$ or
$u_{1x}=u_{0x}=0$. If $u_{1x}=u_{0x}=0$, then it follows from
equation (4.40) that $h_y=0$. Thus, we have $h=0$ because
$h(0,0)=0$. Therefore, $u,v,h$ are linearly dependent.\\
If $d=3$, then comparing the coefficients of $z^2,z$ and $z^0$ of
equation (4.40), we have the following equations:
\begin{equation}
\nonumber
  \left\{ \begin{aligned}
  u_{1x}v_{1y}-3u_3h_x = 0~~~~~~~~~~~~~~~~~~~(4.41) \\
  u_{1x}v_{0y}-v_{1y}u_{0x}-2u_2h_x = 0~~~~~~~(4.42) \\
  u_{0x}v_{0y}-u_1h_x-v_1h_y = 0~~~~~~~~~~~(4.43)\\
                          \end{aligned} \right.
\end{equation}
It follows from equations (4.3) and (4.5) that
$-h_x(v_{1y}z+v_{0y})(3u_3z^2+2u_2z+u_1)+[h_x(u_{3y}z^3+u_{2y}z^2+u_{1y}z+u_{0y})-h_y(u_{1x}z+u_{0x})]v_1=0$
~~~~~~~~~~~~~~~~(4.44) \\
Comparing the coefficients of $z^3$ of the above equation, we have
$h_x(3v_{1y}u_3-u_{3y}v_1)=0$. Thus, we have $h_x=0$ or
$3u_3v_{1y}=v_1u_{3y}$.

(1) If $h_x=0$, then it follows from equation (4.44) that
$h_y(u_{1x}z+u_{0x})=0$. That is, $h_y=0$ or $u_{1x}=u_{0x}=0$. If
$u_{1x}=u_{0x}=0$, then it follows from equation (4.40) that
$h_y=0$. Thus, we have $h=0$ because $h(0,0)=0$. Therefore, $u,v,h$
are linearly dependent.

(2) If $3u_3v_{1y}=v_1u_{3y}$, then
$$\frac{u_{3y}}{u_3}=3\frac{v_{1y}}{v_1}~~~~~~~~~~~~~~~~~~(4.45)$$
If $v_{1y}=0$, then $u_{3y}=0$. It follows from equation (4.41) that
$h_x=0$. Then it follows from the arguments of Case II (1) that
$u,v,h$ are linearly dependent. We can assume that $u_{3y}v_{1y}\neq
0$. Then we have that $u_3=e^{\bar{d}(x)}v_1^3$ by integrating the
two sides of equation (4.45) with respect to $y$, where $\bar{d}(x)$
is a function of $x$. Since $u_3,v_1\in {\bf K}[x,y]$, we have
$e^{\bar{d}(x)}\in {\bf K}[x]$. That is,
$$u_3=d(x)v_1^3~~~~~~~~~~~~~~~~~~~~~~~~~~~(4.46)$$
with $d(x)\in {\bf K}[x]$, $d(x)\neq 0$. Substituting equations
(4.6) and $(4.46)$ to equation (4.41), we have that
$$-3d(x)v_1^3h_x=v_{1y}^2~~~~~~~~~~~~~~~~~(4.47)$$
Then we have $v_{1y}=0$ by comparing the degree of $y$ of equation
(4.47). It follows from equation (4.47) that $h_x=0$. Then it
follows from the arguments of Case II(1) that $u,v,h$ are linearly
dependent.\\

(b) If $v_{1x}=u_{dy}=0$, then it follows from equation (4.6) that
$u_d\in {\bf K}^*$. then equation (4.7) has the following form:
$(u_{1x}z+u_{0x})(v_{1y}z+v_{0y})-v_{0x}(u_{(d-1)y}z^{d-1}+\cdots+u_{1y}z+u_{0y})-h_x(du_dz^{d-1}+\cdots+u_1)-h_yv_1=0$
~~~~~~~~~~~~~~~~~~~~(4.48)\\
Then we have the following equations:
$$-v_{0x}u_{iy}-(i+1)u_{i+1}h_x=0~~~~~~~~~~~~~~~~~~~~(4.49)$$
by comparing the coefficients of $z^i$ of equation (4.49) for
$i=d-1,d-2,\ldots,3$. Comparing the coefficients of $z^2,z$ and
$z^0$ of equation (4.48), we have the following equations:
\begin{equation}
\nonumber
  \left\{ \begin{aligned}
  u_{1x}v_{1y}-v_{0x}u_{2y}-3u_3h_x = 0~~~~~~~~~~~~~~~~~~~(4.50) \\
  u_{1x}v_{0y}+v_{1y}u_{0x}-v_{0x}u_{1y}-2u_2h_x = 0~~~~~~~(4.51) \\
  u_{0x}v_{0y}-v_{0x}u_{0y}-u_1h_x-v_1h_y = 0~~~~~~~~~~~(4.52)\\
                          \end{aligned} \right.
\end{equation}
It follows from equations (4.3) and (4.5) that\\
$[v_{0x}h_y-h_x(v_{1y}z+v_{0y})](du_dz^{d-1}+(d-1)u_{d-1}z^{d-2}+\cdots+u_1)+[h_x(u_{(d-1)y}z^{d-1}+\cdots+u_{1y}z+
u_{0y})-h_y(u_{1x}z+u_{0x})]v_1=0$~~~~~~~~~~~~~~~~~~~~~~~(4.53)\\
Then we have $h_xv_{1y}=0$ by comparing the coefficients of $z^d$ of
equation (4.53). That is, $h_x=0$ or $v_{1y}=0$.

(3) If $h_x=0$, then equation (4.53) has the following form :
$$v_{0x}h_y(du_dz^{d-1}+(d-1)u_{d-1}z^{d-2}+\cdots+u_1)-h_y(u_{1x}z+u_{0x})v_1=0~~~~~~~~~~~~~~(4.54)$$
Comparing the coefficients of $z^{d-1}$ of equation (4.54), we have
that $v_{0x}h_y=0$. That is, $v_{0x}=0$ or $h_y=0$. If $v_{0x}=0$,
then it reduces to Case II(a). If $h_y=0$, then $h=0$ because
$h(0,0)=0$. Thus, $u,v,h$ are linearly dependent.

(4) If $v_{1y}=0$, then comparing the coefficients of $z^{d-1}$ and
$z^0$ of equation (4.53), we have
$$(v_{0x}h_y-h_xv_{0y})du_d+h_xu_{(d-1)y}v_1=0~~~~~~~~~~~~~~~~~(4.55)$$
and
$$(v_{0x}h_y-h_xv_{0y})u_1+(h_xu_{0y}-h_yu_{0x})v_1=0~~~~~~~~~~~(4.56)$$
If $d>3$, then it follows from equation (4.49)$(i=d-1)$ that
$$h_x=-\frac{1}{du_d}v_{0x}u_{(d-1)y}~~~~~~~~~~~~~~~~~(4.57)$$
If $d=3$, then it follows from equation (4.50) that
$$h_x=-\frac{1}{3u_3}v_{0x}u_{2y}~~~~~~~~~~~~~~~~~~~~~(4.57)$$
Substituting equation (4.57) to equation (4.55), we have
$$v_{0x}[h_y-\frac{v_1}{d^2u_d^2}u_{(d-1)y}^2+\frac{1}{du_d}u_{(d-1)y}v_{0y}]=0$$
for $d\geq 3$. Thus, we have $v_{0x}=0$ or
$$h_y=\frac{v_1}{d^2u_d^2}u_{(d-1)y}^2-\frac{1}{du_d}u_{(d-1)y}v_{0y}~~~~~~~~~~~(4.58)$$
If $v_{0x}=0$, then it reduces to Case II(a). Otherwise,
substituting equations (4.57) and (4.58) to equation (4.52),we have
that
$$u_{0x}v_{0y}-v_{0x}u_{0y}=-\frac{u_1}{du_d}v_{0x}u_{(d-1)y}+\frac{v_1^2}{d^2u_d^2}u_{(d-1)y}^2-\frac{v_1}{du_d}u_{(d-1)y}v_{0y}~~~~~~~(4.59)$$
Substituting equations (4.57) and (4.58) to equation (4.56), we have
that
$$\frac{u_1v_1}{d^2u_d^2}u_{(d-1)y}^2v_{0x}-\frac{v_1^2}{d^2u_d^2}u_{(d-1)y}^2u_{0x}+\frac{v_1}{du_d}u_{(d-1)y}(u_{0x}v_{0y}-v_{0x}u_{0y})=0~~~~~(4.60)$$
Then we have $\frac{v_1^3}{d^3u_d^3}u_{(d-1)y}^3=0$ by substituting
equations (4.6) and (4.59) to equation (4.60). That is,
$u_{(d-1)y}=0$. It follows from equation (4.49)$(i=d-1)$ that
$h_x=0$. Then it reduces to Case II (3).\\

(c) If $u_{dy}=u_{(d-1)y}=0$, then it follows from equations (4.2)
and (4.5) that
$(u_{1x}z+u_{0x})(v_{1y}z+v_{0y})-(v_{1x}z+v_{0x})(u_{(d-2)y}z^{d-2}+\cdots+u_{1y}z+u_{0y})-h_x(du_dz^{d-1}+\cdots+u_1)-h_yv_1=0$
~~~~~~~~~~~~~~~~~~~~~~~~~~~~~~~~~~~~~(4.61) \\
Comparing the coefficients of $z^j$ of equation (4.61) for
$j=d-2,\ldots,3$, we have the following equations:
$$-v_{1x}u_{(j-1)y}-v_{0x}u_{jy}-(j+1)u_{j+1}h_x=0~~~~~~~~~~~~~~~~~(4.62)$$
Comparing the coefficients of $z^{d-1},z^2,z$ and $z^0$, we have the
following equations:
\begin{equation}
\nonumber
  \left\{ \begin{aligned}
  -v_{1x}u_{(d-2)y}-du_dh_x=0~~~~~~~~~~~~~~~~~~~~~~~~~~~~~~~~~~~~(4.63)\\
  u_{1x}v_{1y}-v_{1x}u_{1y}-v_{0x}u_{2y}-3u_3h_x = 0~~~~~~~~~~~~~~~~~~~(4.64) \\
  u_{1x}v_{0y}+v_{1y}u_{0x}-v_{1x}u_{0y}-v_{0x}u_{1y}-2u_2h_x = 0~~~~~~~(4.65) \\
  u_{0x}v_{0y}-v_{0x}u_{0y}-u_1h_x-v_1h_y = 0~~~~~~~~~~~~~~~~~~~~~~~(4.66)\\
                          \end{aligned} \right.
\end{equation}
It follows from equations (4.3) and (4.5) that\\
$[(v_{1x}z+v_{0x})h_y-h_x(v_{1y}z+v_{0y})](du_dz^{d-1}+(d-1)u_{d-1}z^{d-2}+\cdots+u_1)+[h_x(u_{(d-2)y}z^{d-2}+\cdots+u_{1y}z+u_{0y})-h_y(u_{1x}z+u_{0x})]v_1=0$
~~~~~~~~~~~~~~~~~~~~~~~~~~~~~(4.67)\\
Comparing the coefficients of $z^d$ and $z^{d-1}$ of equation
(4.67), we have the following equations:
\begin{equation}
\nonumber
  \left\{ \begin{aligned}
  du_d(v_{1x}h_y-h_xv_{1y})=0~~~~~~~~~~~~~~~~~~~~~~~~~~~~~~~~~~~~~~\\
  (d-1)u_{d-1}(v_{1x}h_y-h_xv_{1y})+du_d(v_{0x}h_y-h_xv_{0y}) = 0 \\
                          \end{aligned} \right.
\end{equation}
That is,  $v_{1x}h_y-h_xv_{1y}=0$ and $v_{0x}h_y-h_xv_{0y}=0$. Then
equation (4.67) has the following form :
$$h_x(u_{(d-2)y}z^{d-2}+\cdots+u_{1y}z+u_{0y})-h_y(u_{1x}z+u_{0x})=0~~~~~~~~~~~(4.68)$$
Then we have $h_xu_{ky}=0$ by comparing the coefficients of $z^k$ of
equation (4.68) for $2\leq k\leq d-2$. Thus, we have $h_x=0$ or
$u_{(d-2)y}=\cdots=u_{2y}=0$. If $u_{(d-2)y}=\cdots=u_{2y}=0$, then
it follows from equation (4.63) that $h_x=0$. If $h_x=0$, then it
follows from equation (4.68) that $h_y=0$ or $u_{1x}=u_{0x}=0$. If
$u_{1x}=u_{0x}=0$, then it follows from equation (4.63) and (4.62)
that $v_{1x}=0$ or $u_{(d-2)y}=\cdots=u_{2y}=0$. If $v_{1x}=0$, then
it reduces to Case II(b). If $u_{(d-2)y}=\cdots=u_{2y}=0$, then it
follows from equation (4.64) that $u_{1y}=0$. It follows from
equation (4.65) that $u_{0y}=0$. Then it follows from equation
(4.66) that $h_y=0$. Thus, we have $h=0$ because $h(0,0)=0$.
Therefore, $u,v,h$ are linearly dependent.
\end{proof}

\begin{cor}
Let $H=(u(x,y,z),v(x,y,z),h(x,y))$ be a polynomial map with
$\deg_zv(x,y,z)\leq 1$. Assume that $H(0)=0$ and the components of
$H$ are linearly independent over ${\bf K}$. If $JH$ is nilpotent,
then \\
$(1)$ $H$ has the form of Theorem 2.2 in the case $\deg_zu=0$ and
$(\deg_yu,\deg_yh)\leq 3$.\\
$(2)$ $H$ has the form of Proposition 3.1 in the case $\deg_zv=0$
and $(\deg_xv,\deg_xh)\leq 3$.\\
$(3)$ there exists $T\in \operatorname{GL}_3({\bf K})$ such that
$T^{-1}HT$ has the form of Theorem 2.2 in the case $\deg_zu\cdot
\deg_zv\geq 1$ and $(\deg(u(x,y,0)-\lambda v(x,y,0)),\deg
h(x,y))\leq 3$ for some $\lambda \in {\bf K}$.
\end{cor}
\begin{proof}
Since $u,v,h$ are linearly independent, so it follows from Theorem
4.1 that $\deg_z u\leq 1$. Thus, we have $\deg_zu\leq 1$ and
$\deg_zv\leq 1$. Then the conclusions (1) and (2) follow from
Theorem 2.2 and Proposition 3.1 respectively. The conclusion of (3)
follows from the proof of Theorem 3.10 in \cite{13} and Theorem 2.2.
\end{proof}


\begin{thebibliography}{99}
\bibitem{2} H. Bass, E. Connell and D. Wright, \newblock {\em The Jacobian Conjecture: Reduction of Degree and Formal Expansion of the Inverse}, Bulletin of the American Mathematical Society,\@ 7 (1982), 287-330.
\bibitem{10} M. de Bondt, \newblock {\em Quasi-translations and counterexamples to the homogeneous dependence problem}, Proceedings of the American Mathematical Society 134 (2006) 2849-2856.
\bibitem{8} M. de Bondt and A. van den Essen, \newblock {\em The Jacobian conjecture: linear triangularization for homogeneous polynomial maps in dimension three}, Report 0413, University of Nijmegen, The Netherlands, 2004.
\bibitem{18} Marc Chamberland and Arno van den Essen, \newblock {\em Nilpotent Jacobian in dimension three}, Journal of Pure and Applied Algebra \@ 205 (2006) 146-155.
\bibitem{3} A. van den Essen, Polynomial Automophisms and the Jacobian Conjecture, Vol. \@ 190 in Progress in Mathematics Birkhauser Basel, \@ 2000.
\bibitem{9} A. van den Essen, \newblock {\em Nilpotent Jacobian matrices with independent rows}, Report 9603, University of Nijmegen, The Netherlands, 1996.
\bibitem{7} E. Hubbers, \newblock {\em The Jacobian conjecture: cubic homogeneous maps in dimension four}, Master's Thesis, University of Nijmegen, The Netherlands, 1994.
\bibitem{1} O.H. Keller, \newblock {\em Ganze Cremona-transformationen Monatschr.}, Math. Phys., \@ 47 (1939) pp.\@ 229-306.
\bibitem{4} S.S.S.Wang, \newblock {\em A Jacobian criterion for separability}, Jour. of Algebra,\@ 65 (1980), 453-494.
\bibitem{6} D. Wright, \newblock {\em The Jacobian conjecture: linear triangularization for cubics in dimension three}, Linear and Multilinear Algebra 34(1993) 85-97.
\bibitem{5} A.V.Yagzhev, \newblock {\em On Keller's problem}, Siberian Math. Journal, \@ 21 (1980), 747-754.
\bibitem{13} Dan Yan, Guoping Tang, \newblock {\em Polynomial maps with nilpotent Jacobians in dimension
three}, Linear Algebra and its Applications, \@ 489 (2016), 298-323.
\bibitem{14} Dan Yan, \newblock {\em Polynomial maps with nilpotent Jacobians in dimension three
I},  Linear Algebra and its Applications(submitted).
\end{thebibliography}
\end{document}